\documentclass[10pt]{amsart}
\usepackage{amsmath,amssymb,verbatim}
\usepackage[utf8]{inputenc}
\usepackage{enumerate}
\usepackage{pdflscape}
\usepackage{amsthm}
\usepackage{mathrsfs}
\usepackage{color}
\usepackage[normalem]{ulem}
\usepackage{cancel}
\usepackage{tikz}
\usetikzlibrary{matrix}
\usepackage[all]{xy}
\usepackage{mathtools}

\newtheorem{theorem}{Theorem}[section]

\newtheorem{proposition}[theorem]{Proposition}

\newtheorem{definition}[theorem]{Definition}
\newtheorem{corollary}[theorem]{Corollary}
\newtheorem{remark}[theorem]{Remark}

\newcommand{\Aut}{\mbox{\rm Aut}}
\newcommand{\Inn}{\mbox{\rm Inn}}

\newcommand{\Z}{{\mathbb Z}}
\newcommand{\Q}{{\mathbb Q}}

\newcommand{\C}{{\mathbb C}}

\newcommand{\SL}{{\rm SL}}

\newcommand{\PSL}{{\rm PSL}}

\newcommand{\GEN}[1]{\langle #1 \rangle}

\newcommand{\U}{\mathcal{U}}

\newcommand{\V}{\mathrm{V}}

\newcommand{\pa}[2]{\varepsilon_{#2}(#1)}

\title{Finite Subgroups of Group Rings: A survey}

\author{Leo Margolis and Ángel del Río}
\thanks{This research is partially supported by FWO (Research Foundation of Flanders), the Spanish Government under Grant MTM2016-77445-P with ``Fondos FEDER'' and, by Fundaci\'{o}n S\'{e}neca of Murcia under Grant 19880/GERM/15.}



\begin{document}

\begin{abstract}
In the 1940's Graham Higman initiated the study of finite subgroups of the unit group of an integral group ring. Since then many fascinating aspects of this structure have been discovered. Major questions such as the Isomorphism Problem and the Zassenhaus Conjectures have been settled, leading to many new challenging problems. In this survey we review classical and recent results, sketch methods and list questions relevant for the state of the art. 
\end{abstract}

\maketitle

\section{Introduction}

Group rings first come up as a natural object in the study of  representations of groups as matrices over fields or more generally as endomorphisms of modules. 
They also appear in topology, knot theory and other areas in pure and applied mathematics. For example, many error correcting codes 
can be realized as ideals in group algebras and this algebraic structure has applications on decoding algorithms.

The aim of this survey is to revise the history and state of the art on the study 
of the finite subgroups of units in group rings with special emphasis on integral group rings of finite groups.
For an introduction including proofs of some first results the interested reader might want to consult \cite{Sehgal1993, delRio18}. Other surveys touching on the topics considered here include \cite{SehgalSurvey03,KimmerleDMV}.

Let $G$ be a group, $R$ a ring and denote by $\U(R)$ the unit group of $R$ and by $RG$ the group ring of $G$ with coefficients in $R$. 
The main problem can be stated as follows:

\begin{quote}
 \textbf{Main Problem}: Describe the finite subgroups of $\U(R G)$ and, in particular, its torsion elements.
\end{quote}

This problem, especially in the case of integral group rings of finite groups, has produced a lot of beautiful results which combine group theory, ring theory, number theory, ordinary and modular representation theory and other fields of mathematics. Several answers to the Main Problem have been proposed. The strongest ones, such as the Isomorphism and Normalizer Problems and the Zassenhaus Conjectures, introduced below, are true for large classes of groups, but today we know that they do not hold in general. Other possible answers, as the Kimmerle or Spectrum Problems are still open. We hope that this survey will stimulate research on these and other fascinating questions on group rings. For this purpose we include several open problems and revise the status of some problems given before in \cite{Sehgal1993} in our final remark.    

One of the main motivations for studying finite subgroups of $RG$ in the case where $G$ is finite is the so called  
Isomorphism Problem which asks whether the ring structure of $RG$ determines the group $G$ up to isomorphism, i.e. 

\begin{quote}
\textbf{The Isomorphism Problem}: 
Does the group rings $RG$ and $RH$ being isomorphic imply that so are the groups $G$ and $H$?\vspace{.3cm}

\noindent
\textbf{(ISO)} is the Isomorphism Problem for $R=\Z$ and $G$ finite.
\end{quote}

Observe that the Isomorphism Problem is equivalent to the problem of whether all 
the group bases of $RG$ are isomorphic. A \emph{group basis} is a group of units in $RG$ which is a basis of $RG$ over $R$.
It is easy to find negative solutions to the Isomorphism Problem if the coefficient ring is big, for example, if $G$ and $H$ are finite then $\C G$ and $\C H$ are isomorphic if and only if $G$ and $H$ have the same list of character degrees, with multiplicities. In particular, if $G$ is finite and abelian then $\C G$ and $\C H$ are isomorphic if and only if $G$ and $H$ have the same 
order. As (ISO) was considered a conjecture for a long time it is customary to speak of counterexamples to (ISO).

The  ``smaller'' the coefficient ring is, the harder it is to find a negative solution for the Isomorphism Problem.
This is the moral of the following:

\begin{remark}\label{ISO-Extension}
If there is a ring homomorphism $R\rightarrow S$ then $SG\cong S\otimes_R RG$. 
Thus a negative answer to the Isomorphism Problem for $R$ is also a negative answer for $S$.
In particular, a counterexample for (ISO) is a negative solution for the Isomorphism Problem for all the rings.
\end{remark}

In the same spirit, at least in characteristic zero, the ``smaller'' the ring $R$ is the harder it is to 
construct finite subgroups of $RG$ besides those inside the group $\U(R)G$ of \emph{trivial units}. 
For example, if $G$ is a finite abelian group then all the torsion elements of $\U(\Z G)$ are trivial, i.e. contained in 
$\pm G$. 
This implies that (ISO) has a positive solution for finite abelian groups. 
This is a seminal result from the thesis of Graham Higman \cite{Higman1940}, where the Isomorphism Problem appeared for 
the first time and which raised the interest in the study of units of integral group rings. More than 20 years later Albert Whitcomb proved (ISO) for metabelian groups \cite{Whitcomb}.

The map $\varepsilon: RG \rightarrow R$ associating each element of $RG$ to the sum of its coefficients 
is a ring homomorphism.
This is called the \emph{augmentation map}. It restricts to a group homomorphism $\U(RG)\rightarrow 
\U(R)$ whose kernel is denoted $V(RG)$ and its elements are 
called \emph{normalized units}. 
Clearly $\U(RG)=\U(R)\times V(RG)$, in particular $\U(\Z G)=\pm V(\Z G)$.
It can be easily shown that if $RG$ and $RH$ are isomorphic there is a normalized isomorphism 
$\alpha:RG\rightarrow RH$, i.e. $\varepsilon(\alpha(x))=\varepsilon(x)$ for every $x\in RG$.

Higman's result on torsion units of integral group rings of abelian groups cannot be generalized to non-abelian groups 
because conjugates of trivial units are torsion units which in general are not trivial. 
A natural guess is that all torsion units in the integral group ring of a finite group are of this form, or 
equivalently every normalized torsion unit is conjugate to an element of $G$. 
Higman already observed that $V(\Z S_3)$ contains torsion units which are not conjugate 
in $\U(\Z S_3)$ to trivial units
($S_n$ denotes the symmetric group on $n$ letters).
Since Higman's thesis was not that well known, this was reproven  many years 
later by Ian Hughes and 
Kenneth Pearson. 
They observed however that all the torsion elements of $V(\Z S_3)$ are conjugate to elements of $S_3$ 
in $\Q S_3$ \cite{HughesPearson72}. 
Motivated by this and Higman's result, Hans Zassenhaus conjectured that this holds for all the integral group 
rings of finite groups 
\cite{Zassenhaus}:
  \begin{quote}
  \textbf{The First Zassenhaus Conjecture (ZC1)}:
  If $G$ is a finite group then every normalized torsion unit in $\Z G$ is conjugate in $\Q G$ to an element of 
$G$.
  \end{quote}
Similar conjectures for group bases and, in general for finite subgroups of $\Z G$, are attributed to Zassenhaus \cite{Sehgal84}:
  \begin{quote}
  \textbf{The Second Zassenhaus Conjecture (ZC2)}: 
  If $G$ is a finite group then every group basis of normalized units of $\Z G$ is conjugate in $\Q G$ to $G$.
  \end{quote}
  \begin{quote}
  \textbf{The Third Zassenhaus Conjecture (ZC3)}: 
  If $G$ is a finite group then every finite subgroup of normalized units in $\Z G$ is conjugate in $\Q G$ to a 
subgroup of $G$.
  \end{quote}
Some support for these conjectures came from the following results: If $H$ is a finite subgroup of $V(\Z G)$ then its 
order divides the order of $G$ \cite{ZmudKurennoi1967} and its elements are linearly independent over $\Q$ 
\cite{Higman1940}. The exponents of $G$ and $V(\Z G)$ coincide. 
The last fact is even true replacing $\Z$  by any ring of algebraic integers \cite{CohnLivingstone}. 

The Second Zassenhaus Conjecture is of special relevance because a positive solution for (ZC2) implies a positive 
solution for (ISO). Actually 
  \begin{center}
   (ZC2) \quad $\Leftrightarrow$ \quad (ISO) + (AUT)
  \end{center}
where (AUT) is the following Problem:
  \begin{quote}
  \textbf{The Automorphism Problem (AUT)}: 
  Is every normalized automorphisms of $\Z G$ the composition of the linear extension of an automorphism of $G$ and the 
restriction to $\Z G$ of an inner automorphism of $\Q G$?
  \end{quote}

In the late 1980s counterexamples to the conjectures started appearing.
The first one, by Klaus Wilhelm Roggenkamp and Leonard Lewy Scott \cite{Scott87, Roggenkamp91,Scott92}, was a metabelian negative solution to (AUT) and hence a counterexample to (ZC2) and (ZC3).
Observe that while (ZC2) fails for finite metabelian groups, (ISO) holds for this class by Whitcomb's result mentioned above. 
So in the 1990s there was still some hope that (ISO) may have a positive 
solution in general, as Higman had already stated in his thesis: ``Whether it is possible for two non-isomorphic groups to have isomorphic integral group rings I do not know, but the results of section 5 suggest that it is unlikely'' \cite{Higman1940Thesis}.
However Martin Hertweck found in 1997 two non-isomorphic groups with isomorphic integral group rings 
\cite{HertweckThesis,Hertweck01}. We elaborate on these questions in Section~\ref{GroupBases}. 

So the only of the above mentioned questions open at the end of the 1990s was (ZC1). 
From the 1980s (ZC1) has been proven for many groups including some important classes as nilpotent or metacyclic groups. 
A lot of work was put in trying to prove it for metabelian groups, but 
metabelian counterexamples were discovered by Florian Eisele and Leo Margolis 
in 2017 \cite{EiseleMargolis17}. Details on (ZC1) are provided in Section~\ref{TorsionUnits}.  

Still both (ISO) and the Zassenhaus Conjectures have a positive solution for many 
significant classes of finite groups and many interesting questions on finite subgroups of $\U(\Z G)$ are still open. 
In the introduction we mention just a few of the most relevant. 
It is not known whether (ISO) has a positive solution 
for finite groups of odd order and Hertweck's method can not work in this case. 
A classical open problem is the following special case of the Isomorphism Problem: 

\begin{quote}\textbf{Modular Isomorphism Problem (MIP)}: Let $k$ be a field of characteristic $p$ and $G$ and $H$ finite $p$-groups. Does $kG \cong kH$ imply $G \cong H$?
\end{quote}

Before we give more details on these questions in Section~\ref{GroupBases}, we revise in Section~\ref{General} some methods to attack the problems mentioned above 
as well as known results on (ZC3), the strongest conjecture on the finite subgroups of $\Z G$.

The Zassenhaus Conjectures were possible answers to the Main Problem and in particular (ZC1) was still standing as a possible answer for torsion units until recently. Since all three Zassenhaus Conjectures have been disproved, maybe it is time to reformulate them as the \emph{Zassenhaus Problems}. Another type of answer was proposed by Wolfgang Kimmerle \cite{Ari}:

\begin{quote}\textbf{Kimmerle Problem (KP)}:
Let $G$ be a finite group and $u$ a torsion element in $\V(\Z G)$. Does there exist a group $H$ which contains $G$ such that $u$ is conjugate in $\Q H$ to an element of $G$?
\end{quote}

Observe that while (ZC1) asks if it is enough to enlarge the coefficient ring of $\Z G$ to obtain that all torsion units are trivial up to conjugation, (KP) allows to enlarge also the group basis. 

Recall that the \emph{spectrum} of a group is the set of orders of its torsion elements. A weaker answer to the Main Problem could be provided by solving the following problem: 

\begin{quote}
\textbf{The Spectrum Problem (SP)}:
If $G$ is a finite group do the spectra of $G$ and $V(\Z G)$ coincide?
\end{quote}

See more about these problems in Section~\ref{SectionAltZC1}.

As $V(\Z G)$ and $G$ have the same exponent, at least the orders of the $p$-elements of $V(\Z G)$ and $G$ coincide.
So, it is natural to ask whether the isomorphism classes of the finite $p$-subgroups of $V(\Z G)$ and $G$ are the same. 
It is even an open question whether every (cyclic) finite $p$-subgroup of $V(\Z G)$ is conjugate in $\Q G$ to a 
subgroup of $G$. Section~\ref{pSubgroups} deals with these and other questions on finite $p$-subgroups of $V(\Z G)$.

Techniques and results from modular representation theory have been very useful in the study of the problems mentioned in this article. On the other hand, questions in modular representation theory about the role of the defect group in a block are related to questions on $p$-subgroups in units of group rings. Here rational conjugacy is not as useful and $p$-adic conjugacy is in the focus, as in the $F^*$-Theorem (Theorem~\ref{F*Theorem}) and in Theorem~\ref{HertweckpTorsion}. There is some hope that this kind of results might be applied e.g. to solve the following question of Scott:

\begin{quote} \textbf{Scott's Defect Group Question} \cite{Scott87}: 
Let $\mathbb{Z}_p$ denote the $p$-adic integers, let $G$ be a finite group and let $B$ be a block of the group ring $\Z_p G$. 
Is the defect group of $B$ unique up to conjugation by a unit of $B$ and suitable normalization?
\end{quote}

It is not clear, and should be regarded as part of the problem, what suitable normalization means if the block is not 
the principal block of the group ring. This problem, and indeed Scott's question in its generality, has been solved by Markus 
Linckelmann in case the defect group is cyclic \cite{Linckelmann96}. 

Scott's question is also of interest since, as has been shown by Geoffrey Robinson \cite{Robinson90, Robinson09}, even 
a weak positive answer to it would provide a proof of the $Z_p^*$-Theorem avoiding the Classification of Finite Simple 
Groups. The first proof of the $Z_p^*$-theorem using the CFSG is due to Orest Artemovich \cite{Artemovich}. 
Here the $Z_p^*$-Theorem means an odd analogous of the famous $Z^*$-Theorem of George Glauberman. \\

We use standard notation for cyclic $C_n$ and dihedral group $D_n$ of order $n$; symmetric $S_n$ and alternating group $A_n$ of degree $n$; and linear groups $\SL(n,q)$, $\PSL(n,q)$, etc. For an element $g$ in a group $G$ we denote by $C_G(g)$ the centralizer and by $g^G$ the conjugacy class of $g$ in $G$.

\section{General finite subgroups}\label{General}

Though the group of units of group rings has been studied for about eighty years there are very few classes of group rings for 
which the group of units has been described explicitly. 
For the case of integral group rings, the interested reader can consult the book by Sudarshan Kumar Sehgal \cite{Sehgal1993} 
and by Eric Jespers and Ángel del Río \cite{GRG1,GRG2}. An overview of Higman's thesis \cite{Higman1940Thesis}, the starting point of the area, may be found in \cite{Sandling1981}.
Actually constructing specific units is not that obvious, except for the trivial units of a group 
ring $RG$, i.e. those in $\U(R)G$. 
For example, there is a famous question, studied at least since the 1960s, attributed to Irving Kaplansky \cite{Kaplansky}, though he refers to a question 1.135 by Dmitrij Smirnov and Adalbert Bovdi from the first edition of the Dniester notebook from 1969 (the first edition is entirely included in \cite{Dniester} and the English translation \cite{DniesterEnglish}):

\begin{quote}
\textbf{Kaplansky's Unit Conjecture}: 
If $G$ is a torsion-free group and $R$ is a field then every unit of $RG$ is trivial.
\end{quote}
Kaplansky's Unit Conjecture is still open and few progress has been made on it besides the case of ordered 
groups for which it is easy to verify (in fact the same proof works for unique product groups).
In contrast with this, the only finite groups for which all the units of $\Z G$ are trivial are the abelian groups of 
exponent dividing $4$ or $6$ and the Hamiltonian $2$-groups. Actually these are the only groups for which 
$\U(\Z G)$ is finite \cite{Higman1940Thesis} (see also \cite[Theorem~1.5.6]{GRG1}). 

From now on $R$ is a commutative ring, $G$ is a finite group and we focus on finite subgroups of $\U(RG)$. 
The commutativity of $R$ allows to identify left and right $RG$-modules by 
setting  $rgm=mrg^{-1}$ for $r\in R$, $g\in G$ and $m$ an element in a left or right $RG$-module. 

Higman proved that if $G$ is an abelian finite group then every torsion unit of $\Z G$ is trivial. 
In the 1970s some authors computed $\U(\Z G)$ for some small non-abelian groups $G$. 
For example, Hughes and Pearson computed $\U(\Z S_3)$ \cite{HughesPearson72} and César Polcino Milies computed 
$\U(\Z D_8)$ \cite{PolcinoMilies73}.
As a consequence of these computations it follows that (ZC3) has a positive solution for $S_3$ and $D_8$.
These early results were achieved by very explicit computations. 

A notion to deal with more general classes of groups is the so called \emph{double action module}:

\begin{definition}
Let $H$ be group and let $\alpha:H\rightarrow \U(RG)$ be a group homomorphism. 
Let $(RG)_\alpha$ be the $R(G\times H)$-module whose underlying $R$-module equals $RG$ and the action by 
elements of $G\times H$ is given by:
  $$(g,h)\cdot m = \alpha(h)mg^{-1} \quad (g\in G, h\in H, m\in (RG)_\alpha.$$
\end{definition}

The connection between the Zassenhaus Conjectures and double action modules relies on the following observations.

\begin{proposition}\label{DoubleActionProposition}
Let $\alpha,\beta:H\rightarrow \U(RG)$ be group homomorphisms. 
Then $(RG)_\alpha$ and $(RG)_\beta$ are isomorphic as $R(G\times H)$-modules if and only if there is $u\in \U(RG)$ such 
that $\beta(h)=u^{-1}\alpha(h)u$ for every $h\in H$.
\end{proposition}

\begin{proposition}\label{MRSWGroups}
Let $G$ and $H$ be groups and let $M$ be a left $R(G\times H)$-module. Then $M$ is isomorphic to a double action 
$R(G\times H)$-module if and only if its restriction to $RG$ is isomorphic to the regular right $RG$-module.
\end{proposition}

For our applications $R$ usually is the ring of integers or the field of rationals, and occasionally the ring $\Z_p$ of 
$p$-adic integers. More precisely, consider a finite subgroup $H$ of $V(\Z G)$. Then the embedding $\alpha: H \hookrightarrow 
\U(\Z G)$ defines a double action $\Z(G\times H)$-module but also a double action 
$\Q(G\times H)$-module $(\Q G)_\alpha$. 
By Proposition~\ref{DoubleActionProposition}, to prove that $H$ is conjugate in $\Q G$ to a subgroup of $G$ we need to 
prove that $(\Q G)_\alpha \cong (\Q G)_\beta$ for some group homomorphism $\beta: H \rightarrow G$.
As two $\mathbb{Q}(G \times H)$-modules are isomorphic if and only if they afford the same character,
the following formula is relevant:
  \begin{equation}\label{DoubleActionCharacter}
   \chi_{\alpha}(g,h) = |C_G(g)|\; \varepsilon_g(\alpha(h)) \quad (g\in G, h\in H).
  \end{equation}
Here $\chi_{\alpha}$ denotes the character 
afforded by $(RG)_\alpha$ for 
an arbitrary homomorphism $\alpha:H\rightarrow \U(R G)$ and 
  $$\varepsilon_g:RG \rightarrow R, \quad \sum_{x \in G} a_x x \mapsto \sum_{x\in g^G} a_x.$$
The element $\varepsilon_g(a)$ is called the \emph{partial augmentation} of $a$ at $g$. 
The partial augmentation has an even more practical role in dealing with the Zassenhaus Conjectures via the following:

\begin{proposition}\label{PropMRSW}\cite{MarciniakRitterSehgalWeiss1987},\cite[Lemma 41.4]{Sehgal1993} 
The following are equivalent for a finite subgroup $H$ of $V(\Z G)$:
\begin{enumerate}
 \item $H$ is conjugate in $\Q G$ to a subgroup of $G$. 
 \item There is a group homomorphism $\varphi:H\rightarrow G$ such that for every $h\in H$ and 
$g\in G$ one has $\varepsilon_g(h)\ne 0$ if and only if $g^G=\varphi(h)^G$. 
\end{enumerate}
\end{proposition}

This theorem has been the cornerstone in the study of the Zassenhaus Conjectures. It is the reason why a lot of research has been deployed to study partial augmentations of torsion units of $\Z G$. We collect here some of the most important results in this direction. The first one is also known as the Berman-Higman theorem, named after Higman and Samuil Berman - probably the two earliest researchers in the field.

\begin{proposition}\label{Partial}
Let $G$ be a finite group and let $u$ be a torsion unit of order $n$ in $V(\Z G)$.
\begin{enumerate}
 \item If $n\ne 1$ then $\varepsilon_1(u)=0$ \cite[Proposition~1.5.1]{GRG1}.
 \item If $u$ lies in a group basis then there is an element $g \in G$ such that $\varepsilon_x(u) \neq 0$ if and only if $x \in g^G$. This is an immediate consequence of the class sum correspondence \cite[IV.1 Theorem]{RoggenkampTaylor}.
 \item\label{PartialDivides} If $\varepsilon_g(u)\ne 0$ then the order of $g$ divides $n$ \cite{HertweckBrauer}.
 \item\label{PartialpAdic} 
  If $\varepsilon_g(u)\ne 0$ and the $p$-part of $u$ is conjugate to an element $x$ of $G$ in $\Z_p G$ then $x$ 	
    is conjugate to the $p$-part of $G$ \cite{Hertweck2006}.
 \item\label{PartialSolvable} If $G$ is solvable then $\varepsilon_g(u)\ne 0$ for some element $g$ of order $n$ in $G$ 
\cite{HertweckOrders}.
\end{enumerate}
\end{proposition}

Proposition~\ref{MRSWGroups} also motivated the introduction of an algorithmic method to study finite subgroups $H$ of $V(\Z G)$ using the 
characters of $G$ to obtain restrictions on the partial augmentations of the elements of $H$. 
Each ordinary character $\chi$ of $G$ extends linearly to a map defined on $\C G$, its 
restriction to $H$ is 
the character $\chi_H$ of a $\C H$-module and we have
  $$\chi_H(h) = \sum_{g^G} \varepsilon_g(h) \chi(g), \quad (h \in H)$$
where $\sum_{g^G}$ represents a sum running on representatives of the conjugacy classes of $G$. 
Therefore for each ordinary character $\psi$ of $H$ we have 
  \begin{equation}\label{HeLPSubgroups}
   \frac{1}{|H|} \sum_{h\in H} \sum_{g^G} \varepsilon_g(h) \chi(g) \overline{\psi(h)} = \GEN{\chi_H,\psi}_H \in 
\Z^{\ge 0}.
  \end{equation}
This can be used in combination with Propositions~\ref{MRSWGroups} and \ref{PropMRSW} to prove or disprove the Zassenhaus Conjectures in some cases. 
The information provided by this on partial augmentations is also information about the characters of double action  modules, by \eqref{DoubleActionCharacter}. 
This sometimes helps to construct specific groups of units and eventually counterexamples to the Zassenhaus Conjectures. See Section~\ref{TorsionUnits} for more details.

In case the subgroup is $p$-regular similar formulas are available for $p$-Brauer characters. 
More precisely, if $H$ is a finite subgroup of $\U(\Z G)$ of order coprime with $p$, $\chi$ is a $p$-Brauer character of $G$ and $\psi$ is an ordinary character of $H$ then 
	\begin{equation}\label{HeLPSubgroups-2}
	\frac{1}{|H|} \sum_{h\in H} \sum^{p'}_{g^G} \varepsilon_g(h) \chi(g) \overline{\psi(h)} \in 
	\Z^{\ge 0}
	\end{equation}
where $\sum\limits^{p'}_{g^G}$ represents a sum running on representatives of the conjugacy classes of $p$-regular elements of 
$G$ \cite{HertweckBrauer, MargolisC4C2}. Actually, these are the only partial augmentations relevant for the application of Proposition~\ref{PropMRSW}, because if $h\in H$ and $g$ is $p$-singular then $\varepsilon_g(h)=0$ (see statement \eqref{PartialDivides} of Proposition~\ref{Partial}).

These formulas are the bulk of the method introduced by Indar Singh Luthar and Inder Bir Singh Passi who used it to prove (ZC1) for $A_5$ \cite{LutharPassi1989}. Later it was generalized by Hertweck and used to prove (ZC1) for some small $\operatorname{PSL}(2,q)$ and to give a new short proof for $S_5$ \cite{HertweckBrauer}. It is nowadays known as the HeLP Method. It consists, roughly speaking, in solving formulas \eqref{HeLPSubgroups} and \eqref{HeLPSubgroups-2} for all irreducible $\chi$ and $\psi$ viewing the $\varepsilon_g(h)$ as unknowns and employing additional properties of these integers such as those given in Proposition~\ref{Partial}.
The method has been implemented for the GAP system \cite{GAP, HeLPPaper} for the case where $H$ is cyclic.

The strongest positive results on (ZC3) were achieved by Al Weiss. The first one was proved before by Roggenkamp and Scott for 
the special case of group bases.

\begin{theorem}\label{WeissPResult}\cite{Weiss1988}\cite[Appendix]{Sehgal1993}
Let $R$ be a $p$-adic ring i.e. the integral closure of $\Z_p$ in a finite extension of the field of fractions of $\Z_p$, and let $G$ be a finite $p$-group. 
Then every finite subgroup of $V(RG)$ is conjugate in the units of $RG$ to a subgroup of $G$.
\end{theorem}

This theorem is an application of a deep module-theoretic result of Weiss \cite{Weiss1988} which strongly restricts the possible structure of double action modules of $p$-adic group rings of $p$-groups.
As a consequence of Theorem~\ref{WeissPResult}, (ZC3) holds for $p$-groups. Actually Weiss proved: 

\begin{theorem}\label{WeissNilpotentResult}\cite{Weiss1991}
(ZC3) holds for nilpotent groups.
\end{theorem}

Next theorem collects some other results on (ZC3).

\begin{theorem}\label{ZC3Results}
(ZC3) holds for $G$ in the following cases:
\begin{enumerate}
\item\label{Valenti} $G=C\rtimes A$ with $C$ cyclic, $A$ abelian and $\gcd(|C|,|A|)=1$ \cite{Valenti94}.
\item $G$ is either $S_4$, the binary octahedral group \cite{DokuchaevJuriaans96}, $A_5$, $S_5$ or $\operatorname{SL}(2,5)$ \cite{DokuchaevJuriaansPolcinoMilies97}.
\item All the Sylow subgroups of $G$ are cyclic \cite{JuriaansMilies00}.
\item $|G|=p^2q$ with $p$ and $q$ primes \cite{Liu08}.
\end{enumerate}
\end{theorem}

As it was mentioned in the introduction, the first counterexamples to (ZC2) and (ZC3) were constructed by Roggenkamp and 
Scott as a negative solution to (AUT) \cite{Roggenkamp91,Scott92}. This counterexample was metabelian and supersolvable. 
Using their methods Lee Klingler gave an easier negative solution for such a group of order 2880 
\cite{Klingler1991}. More negative solutions were later constructed by Hertweck \cite{Hertweck02, Hertweck02Illinois}, the 
smallest of order 96 \cite{HertweckHabil}, using groups found by Peter Blanchard as semilocal negative solutions 
\cite{Blanchard}. 

We close this section with a very general problem posed by Kimmerle at a conference \cite{Ari} for which little is known:

\begin{quote}
\textbf{The Subgroup Isomorphism Problem (SIP)}:
What are the finite groups $H$ satisfying the following property for all finite groups $G$?
If $V(\Z G)$ contains a subgroup isomorphic to $H$ then $G$ contains a subgroup isomorphic to 
$H$.
\end{quote}

Note that (SP) is the specification of (SIP) to cyclic groups.
The only groups for which a positive solution for (SIP) has been proven are cyclic 
$p$-groups \cite{CohnLivingstone}, $C_p\times C_p$ for $p$ a prime \cite{KimmerleC2C2,HertweckNonCyclic} and $C_4\times 
C_2$ \cite{MargolisC4C2}. All the known negative solutions to (SIP) are based on Hertweck's counterexample to (ISO).

\section{Group bases}\label{GroupBases}

As already mentioned in the introduction, a lot of research on the units of group rings originally focused on the role 
of group bases inside the unit group. This is directly related to questions such as (ISO) or (ZC2). Still it 
turned out to be very complicated to achieve results for big classes of groups, apart from metabelian groups. 
Roggenkamp and Scott \cite{RoggenkampScott1987} proved (ISO) for finite $p$-groups.
In fact they proved that inside the $p$-adic group ring of a finite $p$-group any two group bases are conjugate and 
hence they are isomorphic. 
This of course implies (ZC2) for this class of groups. 
The stronger results of Weiss \cite{Weiss1988}, quoted in Theorems~\ref{WeissPResult} and 
\ref{WeissNilpotentResult}, were obtained using different methods. 
After these relevant achievements other positive results for (ISO) were obtained by some authors. 
Next theorem summarizes some of the most important classes of solvable groups for which (ISO) has been proved.

\begin{theorem} (ISO) has a positive answer for the following classes of groups:
\begin{enumerate}
\item Abelian-by-nilpotent groups \cite{RoggenkampScott1987}, 
\item Supersolvable groups \cite{KimmerleHabil},
\item Frobenius groups and $2$-Frobenius groups \cite{KimmerleHabil},
\item nilpotent-by-abelian groups (a result of Kimmerle given in \cite[Section XII]{RoggenkampTaylor}).
\end{enumerate}
\end{theorem}

Stronger (yet technical) versions of the first and last statements in the previous theorem can be found in 
\cite{RoggenkampZimmermann92} and \cite{Hertweck92, KimmerleRoggenkampProjective}, respectively.

Another problem about the natural group basis of an integral group ring, which is deeply connected to the solution of (ISO), is the 
so-called Normalizer Problem. Note that the group basis $G$ is obviously normalized by $G$ 
itself and the central units of $\Z G$. The Normalizer Problem asks if these two groups already fill out the normalizer of $G$ in $\U(\Z G)$:

\begin{quote}
\textbf{The Normalizer Problem (NP)}: Let $G$ be a finite group. Is it true that the normalizer of $G$ in the units 
of $\Z G$ is the group generated by $G$ and the central units of $\Z G$?
\end{quote}

For many decades this was called the Normalizer Conjecture and so it is reasonable to speak of counterexamples to (NP). In this survey we concentrate on (NP) for finite soluble groups, see \cite{VanAntwerpen18} for recent results on other classes of groups.
One first important contribution, in a more general context, was given already in the 1960's.

\begin{theorem}\label{Coleman} (Coleman Lemma) \cite{Coleman} 
Let $H$ be a $p$-subgroup of the finite group $G$ and let $R$ be a ring in which $p$ is not invertible. 
Then $N_{\U(RG)}(H) = N_G(H) \cdot C_{\U(R G)}(H)$. In particular, (NP) has a positive solution for 
$p$-groups.
\end{theorem}

A further important contribution by Stefan Jackowski and Zbigniew Marciniak is the following.

\begin{theorem}\label{JackowskiMarciniak}\cite{JackowskiMarciniak87}
Let $G$ be a finite group with normal Sylow $2$-subgroup. Then (NP) has a positive solution for $G$. In particular, it has a positive solution for groups of odd order.
\end{theorem}
 
Finally Hertweck constructed in his thesis \cite{HertweckThesis} counterexamples to (NP) and (ISO).

\begin{theorem}\label{CounterNormalizer}\cite[Theorem A]{Hertweck01}
There is a metabelian counterexample to (NP) of order $2^{25} \cdot 97^2$.
\end{theorem}

\begin{theorem}\cite[Theorem B]{Hertweck01}
There are counterexamples to (ISO) of derived length $4$ and order $2^{21} \cdot 97^{28}$.
\end{theorem}
 
It is not a coincidence that both counterexamples appeared at the same time. 
Actually, to construct his counterexample for (ISO) Hertweck first constructed a 
counterexample $G$ to (NP). This $G$ is different from the group described in Theorem~\ref{CounterNormalizer} and actually is not metabelian. 
He explicitly constructed a unit $t$ in $\Z G$ normalizing $G$ and not acting as an inner 
automorphism of $G$. He then defined an action of an element $c$ on $G$ which is inverting $t$, i.e. $t^c = 
t^{-1}$, and proceeded to show that $X =  G \rtimes \langle c \rangle$ and $Y = \langle G, tc \rangle$ are two 
non-isomorphic group bases of $\Z X$.
In view of this construction and Theorem~\ref{JackowskiMarciniak} it becomes clear that part (a) of the following problem is wide open, since no counterexample to (NP) can serve as a starting point for the construction of 
a counterexample as carried out by Hertweck.

\begin{quote}\textbf{Problem 1:} Does (ISO) have a positive answer for the following classes of groups?
\begin{enumerate}
\item Groups of odd order.
\item Groups of derived length 3.
\end{enumerate}
\end{quote}
 
But even in the case where the order of the group is even (ISO) has ``almost'' a positive answer. 
Namely, any group $G$ can be extended by an elementary abelian group $N$ such that (ISO) has a positive answer for $N 
\rtimes G$. This is a consequence of a strong result obtained by Roggenkamp and Scott: the $F^*$-Theorem. 
See \cite{Hertweck16} for some history of the theorem and also a complete proof of the most general case. 

To state the $F^*$-Theorem let $I_R(G)$ denote the augmentation ideal of a group ring $RG$, i.e. the kernel 
of its augmentation map.
\begin{theorem}\label{F*Theorem}[$F^*$-Theorem]
Let $R$ be a $p$-adic ring and $G$ a finite group with a normal $p$-subgroup $N$ 
containing the centralizer of $N$ in $G$. 
Let $\alpha$ be an automorphism of $RG$ such that $\alpha$ stabilizes $I_R(G)$ and 
$I_R(N)G$. Then $G$ and $\alpha(G)$ are conjugate inside the units of $RG$.
\end{theorem}

If one is only interested in the case of integral coefficients then this can be used to answer (ZC2):

\begin{corollary}\cite[Theorem 1.1]{HertweckKimmerle02}
Let $G$ be a finite group with normal $p$-subgroup $N$ such that the centralizer of $N$ in $G$ is contained in $N$. Then (ZC2) holds for $G$. 
\end{corollary}

Though most of the time the questions mentioned in this section have been studied for special classes of solvable groups, also (almost) simple groups were partly in the focus of attention. We mention some results. The proof of the first theorem uses the Classification of Finite Simple Groups.

\begin{theorem}\cite{KimmerleLyonsSandlingTeague90}
If $\mathbb{Z}G \cong \mathbb{Z}H$ then $G$ and $H$ have isomorphic chief series. In particular, (ISO) holds for finite simple groups. 
\end{theorem}

\begin{theorem}\cite{Petersen76, Bleher95, BleherHissKimmerle95, BleherGeckKimmerle97, Bleher99, BleherKimmerle00}
(ZC2) holds for symmetric groups, minimal simple groups, simple groups of Lie type of small rank, 18 sporadic simple groups and Coxeter groups.
\end{theorem}

In view of these results it might be surprising that all three Zassenhaus Conjectures remain open for alternating groups, cf. \cite[Problem 14]{Sehgal1993}.

We close this section by shortly considering the general Isomorphism Problem.
Sam Perlis and Gordon Walker proved that the Isomorphism Problem for finite abelian groups and rational coefficients 
has a positive solution \cite{PerlisWalker1950}.
Observe that this implies Higman's answer to (ISO) for abelian groups by Remark~\ref{ISO-Extension}.
Richard Brauer asked in \cite{Brauer1963} the following strong version of the Isomorphism Problem: 
Can two non-isomorphic finite groups have isomorphic group algebras over every field? 
Two metabelian finite groups satisfying this were exhibited by Everett Dade 
\cite{Dade1971}. 
This contrasts with the positive result on (ISO) for metabelian groups mentioned above which was already known at the time.
Note that questions on degrees of irreducible complex characters, as 
presented e.g. in \cite[Section 12]{Isaacs1976} or in more recent work on local-global conjectures such as 
\cite{MalleSpaeth16, NavarroMcKay}, can be regarded as questions on what determines the isomorphism type of a complex group algebra. 
Recently a variation of the Isomorphism Problem for twisted group rings has been introduced \cite{MargolisSchnabel}.  

In contrast to the problems described before, (MIP) deals with an object which is finite, but whose unit group fills up almost the whole group algebra. Though extensively studied the problem is only solved when $G$ is either not too far from being abelian or when its order is not too big. Major contributions were given by, among others, Donald Steven Passman, Robert Sandling and Czes{\l}aw Bagi{\'n}ski. We refer to \cite{ HertweckSoriano06, BaginskiKonovalov, EickKonovalov11} for an overview of known results and for a list of invariants of any group basis determined by the modular group algebra. To our knowledge it is not even clear if the choice of the base field $k$ might make a difference for (MIP).

\section{Torsion units - (ZC1)}\label{TorsionUnits}

In Section~\ref{General} we have observed the relevance of partial augmentations for the study of finite subgroups of 
$V(\Z G)$. When studying the First Zassenhaus Conjecture this has even a nicer form: 

\begin{theorem}\cite{MarciniakRitterSehgalWeiss1987}\label{MRSWUnits}
Let $G$ be a finite group and let $u$ be an element of order $n$ in $V(\Z G)$. Then $u$ is conjugate in $\Q G$ to an 
element of $G$ if and only if for every divisor $d$ of $n$ and every $g\in G$ one has $\varepsilon_g(u^d)\ge 0$.
\end{theorem}

Observe that the condition in the last theorem is equivalent to the following: for every $d\mid n$ there is a conjugacy 
class of $G$ containing all the elements at which $u^d$ has non-zero partial augmentation.
 
Most of the early papers on the First Zassenhaus Conjecture dealt with special classes of metacyclic and cyclic-by-abelian groups. For example, (ZC1) was proved for groups of the form $C\rtimes A$ with $C$ and $A$ cyclic of 
coprime order in \cite{PolcinoMiliesSehgal84,PolcinoMiliesRitterSehgal86}. 
This was generalized in \cite{LutharTrama90} for the case where $A$ is abelian (also of order coprime to the order of 
$C$). The proof of the stronger statement in Theorem~\ref{ZC3Results}.\eqref{Valenti} uses these results. 
More positive answers to (ZC1) for special cases of cyclic-by-abelian groups appeared in
\cite{Mitsuda86,MarciniakRitterSehgalWeiss1987,LutharSehgal98,delRioSehgal06}.
Finally Hertweck proved (ZC1) for metacyclic groups in \cite{Hertweck2008}. 
Actually he proved it for groups of the form $G=CA$ with $C$ a cyclic normal subgroup of $G$ and $A$ an abelian subgroup. This was generalized by Mauricio Caicedo and the authors who proved (ZC1) for cyclic-by-abelian groups \cite{CaicedoMargolisdelRio2013}. 
This and Theorem~~\ref{ZC3Results}.\eqref{Valenti} suggest to study the following:
 \begin{quote}
 \textbf{Problem 2:} Does (ZC3) hold for cyclic-by-abelian groups? 
 \end{quote}

Meanwhile (ZC1) was proved for groups of order at most 144, many groups of order less than 288 \cite{HoefertKimmerle06, HermanSingh15, BaechleHermanKonovalovMargolisSingh17} and many other groups. The following list includes the most relevant families of groups for which (ZC1) has been proven:
\begin{itemize}
\item \underline{Metabelian}:
\begin{itemize}
\item $A \rtimes \langle b \rangle$ where $A$ is abelian and $b$ is of prime order smaller than any prime dividing $|A|$ \cite{MarciniakRitterSehgalWeiss1987}. 
\item Groups with a normal abelian subgroup of index $2$ \cite{LutharPassi92}.
\item Cyclic-by-abelian groups \cite{CaicedoMargolisdelRio2013}.
\end{itemize}
\item \underline{Solvable non-metabelian}:
\begin{itemize}
\item Nilpotent groups \cite{Weiss1991}.
\item Frobenius groups of order $p^aq^b$ for $p$ and $q$ primes \cite{JuriaansMilies00}.
\item $P\rtimes A$ with $P$ a $p$-group and $A$ an abelian $p'$-group \cite{Hertweck2006}.
\item $A\times F$ with $A$ abelian and $F$ a Frobenius group with complement of odd order \cite{BaechleKimmerleSerrano18}.
\end{itemize}
\item \underline{Non-solvable}:
\begin{itemize}
\item $A_5$ \cite{LutharPassi1989}, $S_5$ \cite{LutharTramaS5}, $A_6$ \cite{HertweckA6}, $\operatorname{GL}(2,5)$ and the covering group of $S_5$ \cite{BovdiHertweck2008}.
\item $\operatorname{PSL}(2,q)$ for $q \leq 25$ or $q \in \{31,32\}$ \cite{Wagner, Hertweck2006, HertweckBrauer, HertweckA6, Gildea13, KimmerleKonovalov2015, BaechleMargolisEdinb, 4primaryII}. 
\item $\operatorname{PSL}(2,p)$ for $p$ a Fermat or Mersenne prime \cite{FermatMersenne}.
\item $\operatorname{SL}(2,p)$ or $\operatorname{SL}(2,p^2)$ for $p$ prime \cite{delRioSerrano18}. 
\item Finite subgroups of division rings \cite{DokuchaevJuriaans96, DokuchaevJuriaansPolcinoMilies97, BaechleKimmerleMargolisDFG}.
\end{itemize}
\end{itemize}

More positive results for (ZC1) can be found in \cite{AllenHobby80, BhandariLuthar83, RitterSehgal83, SehgalWeiss86, Fernandes87,  BovdiHoefertKimmerle04, MargolisdelRioCW1, MargolisdelRioPAP, MargolisdelRioCW3}. 

As evident from the above, part (a) of the following problem has seen little advances since being included in \cite[Problems 10, 14]{Sehgal1993}.

\begin{quote}\textbf{Problem 3:} Is (ZC1) true for the following groups?
\begin{enumerate}
\item Alternating and symmetric groups.
\item $\operatorname{PSL}(2,p)$ for $p$ a prime.
\end{enumerate}
\end{quote}

As it was mentioned in the introduction a metabelian counterexample to (ZC1) was discovered recently by Eisele and Margolis \cite{EiseleMargolis17}. It is worth to give some explanations on how this counterexample was discovered. Many of the groups for which (ZC1) was proved contained a normal subgroup $N$ such that $N$ and $G/N$ have nice properties (cyclic, abelian or at least nilpotent). Often the proof separates the case where the torsion unit $u$  maps to 1 by the natural homomorphism $\omega_N: \Z G \rightarrow \Z(G/N)$. We write
$$ V(\Z G,N) = \{u \in \V(\mathbb{Z}G) \ | \ \omega_N(u) = 1 \}.$$
The following particular case of (ZC1) was proposed in \cite{Sehgal1993}
\begin{quote}
\textbf{Sehgal's 35th Problem}: 
If $G$ is a finite group and $N$ is a normal nilpotent group of $G$, is every torsion element of $V(\Z G,N)$ conjugate in $\Q G$ to an element of $G$?
\end{quote}
 
The following result of Hertweck, which appeared in \cite{MargolisHertweck}, has interest in itself but it is also important for its applications to Sehgal's 35th Problem, due to statement \eqref{PartialpAdic} of Proposition~\ref{Partial}.

\begin{theorem}\label{HertweckpTorsion}
	Let $G$ be a finite group with normal $p$-subgroup $N$.
	Let $u$ be a torsion element in $\V(\Z G,N)$.
	Then $u$ is conjugate in $\Z_p G$ to an element of $N$.
\end{theorem}

Indeed, it implies that if $u$ is a torsion unit in $V(\Z G,N)$, for $N$ a nilpotent normal subgroup of $G$, then $N$ contains an element $n$ such that for every prime $p$ the $p$-parts of $u$ and $n$ are conjugate in $\Z_p G$. Moreover, by Proposition~\ref{Partial}.\eqref{PartialpAdic}, if $\varepsilon_g(u)\ne 0$ for some $g\in G$ then the $p$-parts of $n$ and $g$ are conjugate in $G$. 

One attempt to attack Sehgal's 35th Problem, already present in \cite{MarciniakRitterSehgalWeiss1987}, is the matrix strategy 
which uses the structure of $\Z G$ as free $\Z N$-module to get a ring homomorphism $\rho:\Z G\rightarrow M_k(\Z N)$, with $k=[G:N]$. Here $M_k$ denotes the $k \times k$-matrix ring.
If $u\in V(\Z G,N)$ then $\rho(u)$ is mapped to the identity via the entrywise application of the augmentation map. 
Using Theorem~\ref{MRSWUnits}, Theorem~\ref{HertweckpTorsion} and a generalization of \eqref{DoubleActionCharacter} it can be proved that if $\rho(u)$ is conjugate in $M_k(\Q N)$ to a diagonal matrix with entries in $N$ then $u$ is conjugate in $\Q G$ to an element of $G$, which would be the desired conclusion. However, Gerald Cliff and Weiss proved that for $N$ nilpotent this approach only works if $N$ has at most one non-cyclic Sylow subgroup \cite{CliffWeiss}. 

Due to this negative result the matrix strategy was abandoned.
However the authors observed in \cite{MargolisdelRioCW1} that some results in the paper of Cliff and Weiss 
can be used to obtain inequalities involving the partial augmentations of torsion elements of $V(\Z G,N)$ which we refer to as the \emph{Cliff-Weiss inequalities}.
In case $N$ is abelian these inequalities take the following friendly form:

\begin{proposition}\cite[Proposition 1.1]{MargolisdelRioCW1}\label{InequalitiesAbelian}
Let $N$ be an nilpotent normal subgroup of $G$ such that $N$ has an abelian Hall $p'$-subgroup $A$ for some prime $p$ and let $u$ be a torsion element of $\V(\Z G, N)$. If $K$ is a subgroup of $A$ such that $A/K$ is cyclic and $n \in N$ then
\[\sum_{g \in nK} |C_G(g)|\pa{u}{g^G} \geq 0. \]
\end{proposition}

The Cliff-Weiss inequalities are actually 
properly stronger than the inequalities \eqref{HeLPSubgroups} for units in $\V(\Z G,N)$ \cite{MargolisdelRioPAP}.
Moreover in \cite{MargolisdelRioCW3} the authors presented an algorithm based on these inequalities and Theorem~\ref{HertweckpTorsion} to search for minimal possible negative solutions to Sehgal's 35th Problem and hence to (ZC1). More precisely the algorithm starts with a 
nilpotent group $N$ and computes a group $G$ containing $N$ as normal subgroup and a list of integers which satisfy the 
Cliff-Weiss inequalities but not the conditions of Theorem~\ref{MRSWUnits}, i.e. they pass the test of the Cliff-Weiss 
inequalities to be the partial augmentations of a negative solution to Sehgal's 35th Problem.

Of course non-trivial solutions of the Cliff-Weiss inequalities do not provide the counterexample yet, because one has to prove the existence of a torsion unit realizing the partial augmentations provided by the algorithm. 
By the double action strategy this reduces to a module theoretical problem, namely one has to prove that there is a certain  $\Z(G\times C_n)$-lattice which is isomorphic to a double action module by Proposition~\ref{MRSWGroups}, where $n$ is the order of the hypothetical unit which is determined by the partial augmentations (see the paragraph after Theorem~\ref{HertweckpTorsion}). A first step to obtain this lattice 
consists in showing the existence of a $\Z_p(G\times C_n)$-lattice with the same character as the double action $\Q(G\times C_n)$-module, which exists since the partial augmentations of the hypothetical unit satisfies the constraints of the HeLP-method, for every prime $p$. 
By the results of Cliff and Weiss, a unit satisfying also the Cliff-Weiss inequalities corresponds to a 
$\Z_p(G\times C_n)$-lattice which is free as $\Z_p N$-lattice. The fundamental ingredient which allows the construction to work at this point is that the $p$-Sylow subgroup $N_p$ of $N$ is a direct factor in $N$. Hence $\Z_p N = \Z_p N_{p'} \otimes_{\Z_p} \Z_p N_p$ and the representation theory of the first factor is easy to control. It turns out that assuming $N$ is abelian and that a 
$\Z_p(G\times C_n)$-lattice 
which is free of rank $1$ as $\Z_p G$-lattice (compare with Proposition~\ref{MRSWGroups}), assuming only 
that the action of $G$ on $N$ satisfies a certain, relatively weak, condition \cite[Section 5]{EiseleMargolis17}.

Once such a $\Z_p(G\times C_n)$-lattice $M_p$ is constructed for every prime $p$, one obtains a $Z_{(\pi)}(G\times C_n)$-lattice with the same character as each $M_p$, where $\Z_{(\pi)}$ denotes the localization of $\Z$ at the set of prime divisors of the order of $G$. So one obtains what is usually called a	semilocal counterexample. It remains to show how this lattice can be ``deformed'' into a $\Z(G\times C_n)$-lattice with the same character. This is done in \cite[Section 6]{EiseleMargolis17} in a rather general context which could be applied also to non-cyclic groups and other coefficient rings. In the situation of (ZC1) this boils down to checking that $G$ does not map surjectively  onto certain groups (which is, in this case, equivalent to the Eichler condition for $\Z G$) and that $D(u)$ has an eigenvalue $1$ for any irreducible $\Q$-representation $D$ of $G$.

With all this machinery set up, to find a counterexample to (ZC1) remains a matter of calculations and it turns out that the candidates constructed as minimal possible negative solutions to Sehgal's 35th Problem in \cite{MargolisdelRioCW3} are in fact negative solutions and as such counterexamples to (ZC1). 

The construction gives rise to the following problem.
\begin{quote}\textbf{Problem 4:} Classify those nilpotent groups $N$ such that Sehgal's 35th Problem has a positive solution for any group $G$ containing $N$ as normal subgroup. 
\end{quote}

By Theorem~\ref{HertweckpTorsion} and \cite{CliffWeiss} the class of groups described in this problem contains those nilpotent groups which have at most one non-cyclic Sylow subgroup, cf. \cite{MargolisdelRioCW1} for details. More technical results for the problem can be found in \cite{MargolisdelRioCW1, MargolisdelRioCW3, MargolisdelRioPAP}. By the counterexamples to (ZC1) there are infinitely many pairs of different primes $p$ and $q$ such that the direct product of a cyclic group of order $p \cdot q$ with itself is not contained in this class. This is particularly the case for $(p,q) = (7,19)$, but not for $(p,q)$ with $p\le 5$. 

The evidence provided by positive solutions to (ZC1) and Problem 3 and by the counterexamples to (ZC1) suggests that the following might have a positive answer: 
\begin{quote}\textbf{Problem 5:}
Is (ZC1) true for supersolvable groups? Is it true at least for supersolvable metabelian groups?
\end{quote}

\section{Torsion units - Alternatives to (ZC1)}\label{SectionAltZC1}

Already when (ZC1) was still open several weaker forms of the conjecture were proposed, and all of these remain open, in a strong sense. The first and strongest is the Kimmerle Problem mentioned in the introduction. It was posed by Kimmerle at a conference in Oberwolfach \cite{Ari}[Problem 22], partly motivated by an observation of Hertweck on the group rings of $\operatorname{PSL}(2,q)$ and $\operatorname{PGL}(2,q)$ \cite[Remark 6.2]{HertweckBrauer}. 
As (ZC1) was regarded as the main problem, (KP) was not studied much by itself. The authors observed in \cite{MargolisdelRioPAP} that (KP) is actually equivalent to Problem 44 from \cite{Sehgal1993}, a generalization of a question posed by Bovdi (see \cite[p. 26]{BovdiBook} or \cite[Problem 1.5]{ArtamonovBovdi1989}). More precisely: 

\begin{proposition}\label{KPGenBP}
Let $G$ be a finite group and $u$ a torsion element of $\V(\Z G)$. Consider $G$ as a subgroup of the symmetric group $S_G$ in the standard way. The following are equivalent:
\begin{enumerate}
\item $u$ is conjugate to an element of $G$ in the rational group algebra of some group containing $G$.
\item $u$ is conjugate to an element of $G$ in the rational group algebra of $S_G$.
\item For every positive integer $m$ different from the order of $u$ the coefficients of $u$ corresponding to elements of $G$ of order $m$ sum up to $0$.
\end{enumerate}
\end{proposition}

We summarize some results on (KP). The first two follow from results of Stanley Orlando Juriaans, Michael Dokuchaev and Sehgal using Proposition~\ref{KPGenBP}.

\begin{theorem}\label{KP} (KP) has a positive answer if one of the following holds.
\begin{enumerate}
\item $G$ is metabelian \cite{DokuchaevSehgal1994},
\item $G$ is solvable, has only abelian Sylow subgroups and $u$ is of prime power order \cite{Juriaans1994},
\item $u$ is of prime order \cite{KimmerleMargolis17},
\item $G$ has a Sylow tower \cite{BaechleKimmerleSerrano18}. In particular, for $G$ supersolvable.
\end{enumerate}
\end{theorem} 

It was observed in \cite{BaechleKimmerleSerrano18} that the counterexamples to (ZC1) constructed in \cite{EiseleMargolis17} can not provide negative solutions to (KP) as they have Sylow towers and are also metabelian.
Actually, as explained above, the methods in \cite{EiseleMargolis17} can probably allow to construct more counterexamples $G$, some of which might not have a Sylow tower and not be metabelian. However all units providing counterexamples with this method will live in $\V(\Z G, N)$ for a normal nilpotent subgroup $N$ of $G$. On the other hand, for elements in $\V(\Z G, N)$, the Kimmerle Problem has a positive solution by Proposition~\ref{KP} and \cite[Theorem 3.3]{MargolisdelRioPAP}. 
In this sense a solution to (KP) would need significant new ideas.

The other weaker version of (ZC1) mentioned in the introduction, i.e. (SP), is very natural and more in the style of questions already asked by Higman.
As a consequence of Proposition~\ref{Partial}.\eqref{PartialSolvable} we have:

\begin{corollary}
(SP) has a positive solution for solvable groups.
\end{corollary}

So (SP) has a positive answer for a very big class of groups, a class for which probably there will never be an argument or algorithm that can tell if a specific group in this class satisfies (ZC1) or not. It is very interesting what this class can give for (KP):

\begin{quote}\textbf{Problem 6:} 
Does (KP) hold for solvable groups?
\end{quote}

One weaker version of (SP) which found some attention was also formulated by Kimmerle \cite{Kimmerle2006}. Recall that the \textit{prime graph}, also called the Gruenberg-Kegel graph, of a group $G$ is an undirected graph whose vertices are the primes appearing as order of elements in $G$ and the vertices $p$ and $q$ are connected by an edge if and only if $G$ contains an element of order $pq$.

\begin{quote}\textbf{The Prime Graph Question (PQ)}:
Let $G$ be a finite group. Do $G$ and $\V(\Z G)$ have the same prime graph?
\end{quote}

The structural advantage is that for (PQ) there is a reduction theorem, 
while this is not the case for any of the other questions given above. Recall that a group $G$ is called \emph{almost simple} 
if there is a non-abelian simple group $S$ such that $G$ is isomorphic to a subgroup of $\Aut(S)$ containing $\Inn(S)$,
and in this case $S$ is called the socle of $G$.

\begin{theorem} \cite{KimmerleKonovalov2016}
Let $G$ be a finite group. Then (PQ) has a positive answer for $G$ if and only if it has a positive answer for all almost simple homomorphic images of $G$.
\end{theorem}

So one might hope that the Classification of Finite Simple Groups can provide a way to prove (PQ) for all groups. But a lot remains to be done, since many series of almost simple groups still need to be handled. We summarize some important results.

\begin{theorem} 
Let $G$ be a finite almost-simple group. Then (PQ) has a positive answer for $G$ if the socle of $G$ is one of the following:
\begin{enumerate}
\item The alternating group $A_n$ \cite{LutharPassi1989, HertweckA6, Salim2011, Salim2013, BaechleMargolisEdinb, BaechleCaicedo17, BaechleMargolisAn}.
\item $\operatorname{PSL}(2,p)$ or $\operatorname{PSL}(2,p^2)$ for some prime $p$ \cite{HertweckBrauer, 4primaryI}.
\item One of $18$ sporadic simple groups \cite{BovdiKonovalovM12, BovdiKonovalovM11, BovdiKonovalovMcL,BovdiKonovalovM22,BovdiKonovalovSuz,BovdiKonovalovM23,BovdiKonovalovHe,BovdiKonovalovRu,BovdiKonovalovHS,
BovdiKonovalovJanko,BovdiKonovalovConway,BovdiKonovalovM24, KimmerleKonovalov2015, MargolisConway, BaechleMargolisAn}.
\item A group whose order is divisible by at most three pairwise different primes \cite{KimmerleKonovalov2016, BaechleMargolisEdinb} or one of many groups whose order is divisible by four pairwise different primes \cite{4primaryII}.
\end{enumerate} 
\end{theorem}

\section{$p$-subgroups}\label{pSubgroups}

In this section we revise the main results and questions on the finite $p$-subgroups of $\U(\Z G)$ for $G$ a finite group and $p$ a prime integer.  
The questions are the specialization to $p$-subgroups of $V(\Z G)$ of the problems given above which we 
refer to by adding the prefix ``$p$-''. 
For example, the $p$-versions of (ZC3) and (SIP) are as follows: 
\begin{quote}
  \textbf{($p$-ZC3)}: 
  Given a finite group $G$, is every finite $p$-subgroup of $V(\Z G)$ conjugate in $\Q G$ to a subgroup of $G$?
\end{quote}
\begin{quote}
  \textbf{($p$-SIP)}.
  What are the finite $p$-groups $P$ satisfying the following property for all finite groups $G$?
  If $V(\Z G)$ contains a subgroup isomorphic to $P$ then $G$ contains a subgroup isomorphic to $P$.
\end{quote}
The following terminology was introduced in \cite{MargolisPSL, KimmerleMargolis17}. One says that $G$ satisfies a \emph{Weak Sylow Like Theorem} when every finite $p$-subgroup of $V(\Z G)$ is isomorphic 
to a subgroup of $G$.

That the role of $p$-subgroups in $\V(\Z G)$ is very special is expressed already by the Lemma of Coleman (Theorem~\ref{Coleman}), which implies a positive solution for ($p$-NP). Also, the result of James Cohn and Donald Livingstone on the exponent of $\V(\Z G)$, mentioned above, is equivalent to a positive solution to ($p$-SP).

By Theorem~\ref{HertweckpTorsion}, the $p$-version of Sehgal's 35th Problem has a positive answer in general. This was in fact already observed earlier by Hertweck \cite{Hertweck2006}.
Moreover, as a consequence of Theorem~\ref{WeissPResult}, (ZC3) holds for $p$-groups and hence ($p$-ZC2) and ($p$-ISO) hold. These latter also follows from the following result:

\begin{theorem}\cite{KimmerleRoggenkamp93}
If $G$ is a finite solvable group and $P$ is a $p$-subgroup of a group basis of $\Z G$ then $P$ is conjugate in $\Q G$ to a 
subgroup of $G$.
\end{theorem}

In the situation of general $p$-subgroups the knowledge is much more sparse. 
There is no counterexample to ($p$-ZC3). Neither is there a general answer to ($p$-ZC1), not even for units of order $p$, though in this case ($p$-KP) holds, cf. Theorem~\ref{KP}.
Note, that all positive results for (SIP), mentioned in Section~\ref{General}, are in fact results for ($p$-SIP). A big step in the solution of this problem might be an answer to the following:

\begin{quote}\textbf{Problem 7:}
Is (SIP) true for elementary-abelian groups?
\end{quote}

We collect here some results on ($p$-ZC3) and Weak Sylow Like Theorems:

\begin{theorem}
If $G$ is a finite group and $p$ a prime integer then ($p$-ZC3) has a positive answer for $G$ in the following cases:
\begin{enumerate}
 \item $G$ is nilpotent-by-nilpotent or supersolvable \cite{DokuchaevJuriaans96}.
 \item $G$ has a normal Sylow $p$-subgroup \cite[41.12]{Sehgal1993}. 
 
 \item $G$ is solvable and the Sylow $p$-subgroups of $G$ are abelian \cite[Proposition~2.11]{DokuchaevJuriaans96}.
 
 \item $G$ is solvable and every Sylow subgroup of $G$ is either abelian or a quaternion $2$-group 
  \cite[Theorem~5.1]{DokuchaevJuriaans96}.
  \item $G$ is a Frobenius group \cite{DokuchaevJuriaans96, DokuchaevJuriaansPolcinoMilies97, JuriaansMilies00, BovdiHertweck2008, KimmerleMargolis17}.
  \item $p=2$, the Sylow $2$-subgroup of $G$ has at most $8$ elements and $G$ is not isomorphic to $A_7$ \cite{BaechleKimmerle11, MargolisC4C2}. 
\end{enumerate}
\end{theorem}

\begin{theorem}
$G$ satisfies a Weak Sylow Like Theorem for $p$-subgroups in the following cases:
\begin{enumerate}
 \item $p=2$ and the Sylow $2$-subgroups of $G$ are either abelian, quaternion \cite{KimmerleSylowLike} or dihedral  \cite{MargolisC4C2}.
 \item $G$ has cyclic Sylow $p$-subgroups \cite{KimmerleC2C2, HertweckNonCyclic}.
\end{enumerate}
\end{theorem}

For $G = \operatorname{PSL}(2,r^f)$ the $p$-subgroups of $\V(\Z G)$ found some attention starting with \cite{HertweckHoefertKimmerle}. It is known today that ($p$-ZC3) has a positive answer for $G$ if $p \neq r$, or $p=r=2$, or $f=1$ \cite{MargolisPSL}. Also a Weak Sylow Like Theorem holds for $G$ if $f \leq 3$ \cite{HertweckHoefertKimmerle, BaechleMargolis15}. \\

\textbf{Remark:}
A quarter of a century ago Sehgal included a list of 56 open problems in his book on units of integral group rings \cite{Sehgal1993}. Several of those concerned topics mentioned in this article. Some have been solved, while others remain open. 
\begin{enumerate}
\item Problem 8 is (SP), Problem 45 is (MIP) and hence both are still open. Problem 43 is (NP) and has been solved by Hertweck.
\item By Proposition~\ref{KPGenBP}, Problem 9 is equivalent to ($p$-KP) and Problem 44 to (KP). Both remain open.
\item Problems 10 and 11 ask to prove (ZC1) and (ZC3) for symmetric groups. Problem 14 asks about (ZC1)-(ZC3) for alternating groups. They remain open. 
\item Problems 33, 34 and 37 are related with the matrix strategy mentioned in Section~\ref{TorsionUnits}. They have been solved by Cliff and Weiss \cite{CliffWeiss}.
\item Problems 12, 35 and 36 are solved negatively by the counterexample to (ZC1) \cite{EiseleMargolis17}.
\item Problem 32 asks the following: Let $N$ be a normal subgroup in $G$, where $G$ is a finite group, and $u$ a torsion element in $\V(\Z G, N)$. Does the order of $u$ necessarily divide the order of $N$? This has positive solution if $N$ is solvable, as we quickly show: Assume that $u$ is of order $n$ and let $N'$ be a minimal normal subgroup of $G$ contained in $N$.  Then $N'$ must be an elementary-abelian $p$-group, for some prime $p$, as $N$ is solvable. If the projection $u'$ of $u$ in $\Z (G/N')$ has order smaller than $u$, then the order of $u'$ is $n/p$ by \cite[Proposition 4.2]{Hertweck2006}. So arguing by induction on the order of $N$ we can answer the problem for solvable $N$.  

We are thankful to the referee who provided the idea for this proof.
\end{enumerate}

\bibliographystyle{plain}
\bibliography{AGTA}

\noindent
adelrio@um.es\newline
Departamento de Matemáticas, Universidad de Murcia, 30100 Murcia, Spain.\newline

\noindent
leo.margolis@vub.be,\newline
Vrije Universiteit Brussel, Department of Mathematics, Pleinlaan 2, 1050 Brussel, Belgium.

\end{document}